\documentclass[a4paper,oneside]{amsart}

\usepackage{lucacommands}
\usepackage[nonewtxmathopt]{newtxmath}  
\usepackage{stackengine}
\usepackage{mdframed}
\begin{document}

\title{An axiomatization of truth and paradoxicality}


\author{Luca Castaldo}
\address{}
\curraddr{}
\email{castaldluca@gmail.com}
\thanks{}




\date{}

\dedicatory{}

\begin{abstract}
This short note introduces a formal system of truth and paradoxicality, outlining the main motivation, and proving its $\omega$-consistency. The system is called \tp, for \emph{Truth and Paradoxicality}.\footnotemark
\end{abstract}
\footnotetext{
	A preliminary version of \tp was presented at the \emph{XII Workshop on Philosophical Logic} (Buenos Aires, August 2023) and at \emph{Vagueness, Truth, and Semantic Indeterminacy} workshop (Turin, November 2023). This note outlines the basic framework. A more detailed presentation, including e.g. a proof-theoretic analysis of \tp and comparisons with related systems (such as Field's INT \cite{fiePower} and Fujimoto \& Halbach's CD \cite{fhClassical}), will be developed in a forthcoming paper.
}

\maketitle

\allowdisplaybreaks
\section{Introduction}
\label{sec:Intro}
	Semantic notions lie at the core of numerous philosophical debates. Yet, formal theories of semantic notions are notoriously constrained by paradoxes. For instance, the well-known Liar paradox reveals an inconsistency between the laws of classical logic and Tarski's schema, or just T-Schema, $\vphi\lra\T\vphi$. The paradigmatic way to show the inconsistency is to consider a liar sentence $\lambda$, obtained by means of some diagonal construction, such that $\lambda\lra\neg\T\lambda$. The instance of the T-schema for $\lambda$ yields an inconsistency.
	
	Due to the existence of paradoxical sentences, truth is often studied alongside other notions that express when certain sentences are \textit{healthy} \cite{bacCan}, not paradoxical. Healthy sentences are those like ``Snow is black'' or ``$0 = 0$'', that can satisfy the T-Schema without leading to contradiction. Notable examples of healthy properties that have been studied and axiomatized simultaneously with truth include \textit{groundedness} \cite{leiDepends, schAxioms, rosShould}, \textit{well-foundedness} \cite{picAlethic, picReferenceAnd}, \textit{determinateness} \cite{fefAxioms, fhClassical, cnClassical}, \textit{significance} \cite{reiRemarks85,reiRemarks} \textit{strong classicality} \cite{fiePower}.

One of the motivating reasons for having a theory of truth along with an healthy notion is precisely the possibility of having a classification of sentences which are not paradoxical, and that can thus satisfy the T-Schema. For example, Andrew Bacon describes ``the project of diagnosing the paradoxes'' as follows:
	\begin{quote}
		The hallmark of this kind of project is to identify some feature common to all problematic instances of T[-Schema] and to accordingly diagnose the potential for T[-Schema] to fail as being due to the presence of this feature. Once this feature is identified, we are in a position to start explaining why instances of T[-Schema] that have the feature are liable to lead to inconsistency, while instances that do not are not. \hfill \cite[p.307]{bacCan}
	\end{quote}
Surprisingly, however, despite the project described by Bacon aims at identifying some feature common to paradoxical sentences, no axiomatic theory of truth and \textit{paradoxicality} has been developed to date.

The notion of paradoxicality has received attention in the recent literature. Yet, it has proven to be somewhat reluctant to formal treatment. In particular, the \textit{naive conception of paradoxicality} has been shown to be prone to paradox. Informally, the naive conception of paradoxicality is based on a very simple and natural thought, namely that one can identify some principle $\Psi$, such that a sentence $\vphi$ is paradoxical just in case the $\vphi$-instance of $\Psi$ leads to a contradiction.\footnote{
	This and other conceptions of semantic paradoxicality are discussed more thoroughly in \cite{rgConceptions}.
	} 
For instance, liar sentences are paradoxical instances of Tarski's schema. 

The issue with this naive conception is that it generates familiar revenge paradoxes: if one tries, in a theory formulated in classical logic, to formalize the claim that a sentence is paradoxical just in case the instance of Tarski's schema for it leads to a contradiction, one obtains an inconsistency; similarly, if one tries, in a theory formulated in a paracomplete logic like Strong Kleene, to formalize the claim that a sentence is paradoxical just in case the instance of the law of excluded middle for it leads to a contradiction, one also obtains an inconsistency.\footnote{
	See, e.g., \cite{mrGeneralized}, \cite{rgParadoxicality}, \cite{rosParadoxicality}, \cite{isReasoning}.
}
Recently, \cite{crAsymmetries} have investigated possible combinations of rules for naive paradoxicality in both a classical and a non-classical setting, establishing various (im)possibility results.

These negative results raise the challenge of refining our naive notion of paradoxicality, and the current paper tries to move the first steps into this direction. It will develop a theory of truth and non-naive paradoxicality, thereby identifying some feature common to all problematic instances of T-Schema.

The theory introduced below builds on ideas from \cite{casFixed}. It is motivated by the following intuitions.
	\begin{itemize}
		\item Truth is thoroughly \textit{compositional}, in the sense that it commutes with all connectives and quantifiers; truth is \textit{transparent}, in the sense that any sentence \vphi is intersubstitutable with ``\vphi is true'' in every context; truth is \textit{consistent}, in the sense that no sentence can be both true and false. In other words, truth is a transparent notion which behaves according to the Strong Kleene (\SK) semantic conditions.
		\item Truth and paradoxicality are disjoint: no paradoxical sentence can be true or false.
		\item Paradoxicality is \textit{quasi-compositional}. That is, there are \textit{base paradoxical sentences}, whose paradoxicality cannot be further analyzed, it's atomic. A typical example of base paradoxical sentence is a liar sentences $\lambda$. Additionally, there are some sentences, such as $\lambda\land 0=0$, whose paradoxicality is \textit{grounded} in its components.   
		\item Like truth, paradoxicality behaves according to \SK, in the following sense. If a sentence contains enough classical information to be declared true or false, then it will be either true or false, and hence not paradoxical. E.g., according to \SK a conjunction is false iff one conjunct is false. Thus, a conjunction $\vphi\land\psi$ with a false conjunct is false, no matter whether the other conjunct is true, false, paradoxical, or neither. 
		\item A sentence is paradoxical iff its negation is. Intuitively, paradoxical sentences cannot be taken to be true or false without yielding a contradiction. But then, under the natural reading of falsity as true negation, and assuming that $\neg\neg\vphi$ is equivalent to \vphi, to say that \vphi can't be true (false) amounts to saying that that $\neg\vphi$ can't be false (true). Thus a sentence is paradoxical iff its negation is.
	\end{itemize}

\section{Preliminaries}
	The language $\lnat$ denotes the language of arithmetic, and the language $\Lmc$ is defined as $\lnat\cup\{\T,\P\}$. For the purposes of this paper, it is convenient to work with a Tait-style language, in which negation is defined in the usual way. In particular, formulae are built up from \Lmc-literals, i.e. $t = s$, $\T t$, $\P t$, $t \neq s$, $\neg\T t$, and $\P t$, by means of $\lor,\land,\exists,\forall$. Negation is extended to all formulas by De Morgan dualities, with the stipulation that $\neg\neg\vphi:=\vphi$. 
	Greek letters $\vphi,\psi,\xi,\ldots$ range over formulae. By an $\Lmc$-expression, we mean a term or a formula of $\Lmc$. The numeral corresponding to the number $n\in\omega$ is denoted by $\ovl n$. We fix a canonical G\"odel numbering of $\Lmc$-expressions and we take Peano arithmetic (\PA) as base syntax theory -- although of course much weaker systems would suffice -- in which one can carry out a primitive recursive (p.r.) formalization of syntactic notions and operations. If $e$ is an $\Lmc$-expression, the G\"odel number of $e$ is denoted by $\#e$ and $\cor{e}$ is the term representing $\#e$ in $\Lmc$. The sets of terms, closed terms, variables, formulae with $n$ free variables, and sentences of $\Lmc$ are elementary and can be represented in $\Lmc$. In practice, we take the following $\lnat$-predicates to abbreviate the equations for the characteristic functions for such sets, respectively: $\mrm{Tm}, \mrm{ClTm}, \mrm{Var}, \mrm{Fml}^n, \mrm{Sent}$. We also include a function symbol $\mrm{num}(x)$ for the standard numeral function, sending a number to the code of its numeral. Moreover, we employ a functional notation $\mrm{val}(x)$ abbreviating the formula representing in $\lnat$ the evaluation function for closed terms.
	
	For simplicity, we extend \Lmc with finitely many function symbols for elementary syntactic operations, including the usual ones on G\"odel numbers:
	\begin{itemize}
 		\item[\ud\star] \ mapping $(\#t, \#s)$ to $\#(t \star s)$, for $\star \in\{=,\neq\}$;
 		\item[\ud*] \ mapping $\#t$ to $\#{*}{t}$, for $* \in\{\T,\P,\neg\T,\neg\P\}$;
 		\item[\ud\circ] \ mapping $(\#\vphi,\#\psi)$ to $\#(\vphi\circ\psi)$, for $\circ \in\{\land,\lor\}$;
 		\item[\ud Q] \ mapping $(\#v,\#\vphi)$ to $\#Q v\vphi$, for $Q \in\{\forall,\exists\}$.
 	\end{itemize}
A function $\ng{}$ mapping (the code of) an arbitrary \vphi to (the code of the definition of) $\neg\vphi$ can be defined in the obvious way.

Let $\mrm{subst} ( x, y, z )$ arithmetically represent the syntactic substitution of a term coded by $y$ for a variable coded by $z$ in an expression coded by $x$; for example, for an expression $e$, $\mrm{sub} (\cor{e}, \cor{t}, \cor{v_k} )$ is the code of $e[t/v_k]$. For readability, we denote $\mrm{subst} ( x, y, z )$ by $x ( y / z )$. For $x \in \mrm{Fml}^1$, $x ( y )$ is a code of the sentence obtained by substituting the numeral for $y$ for the unique free variable in the formula coded by $x$, that is, $x ( \num y / z )$ for the unique $z$ with $\mathrm{Free} ( z, x )$, where $\mrm{Free} ( z,x )$ arithmetically represents the relation that holds between a code $z$ of a free variable in a formula coded by $x$. Moreover, for a formula $\vphi ( v )$, we define $\cor{\varphi(\dot{x})}:=\cor{\vphi(v)}( \num x / \cor v)$. Moreover, for a formula $\vphi ( v )$, we define $\cor{\vphi(\dot{x})}:=\cor{\vphi(v)}( \mrm{num}(x) / \cor v)$.

For the sake of readability, we write $\T\vphi$ and $\P\vphi$ instead of $\T\cor\vphi$ and $\P\cor\vphi$, respectively. Similarly, we write $\T\vphi(x)$ and $\P\vphi(x)$ instead of $\T\cor{\vphi(\dot{x})}$ and $\P\cor{\vphi(\dot{x})}$, respectively.

We work with a two-sided sequent calculus. A sequent is an expression of the form \GD, for $\Gamma$ and $\Delta$ finite sets of sentences.\footnote{
	For details on sequent calculi, see, e.g., \cite{tsBasic}.
} 
The expression $\Gamma\Leftrightarrow\Delta$ is used as shorthand for the two sequents $\Gamma\Ra\Delta$ and $\Delta\Ra\Gamma$. A double line between two sequents, as in
			\begin{IEEEeqnarray*}{L}
	        \def\fCenter{ \mbox{ $\Ra$ } }
                \ax{ \Gamma\fCenter\Delta }
                \dl
                    \un{\Gamma'\fCenter\Delta'}
		      \Dis
	       \end{IEEEeqnarray*}
indicates that the lower sequent is derivable from the upper sequent via a series of inferences. The notation
			\begin{IEEEeqnarray*}{L}
	        \def\fCenter{ \mbox{ $\Ra$ } }
                \ax{ \Gamma\fCenter\Delta }
                \dl
                \rl{$\mathit{I}$}
                    \un{\Gamma'\fCenter\Delta'}
		      \Dis
	       \end{IEEEeqnarray*}
indicates that the lower sequent is derivable from the upper sequent via an application of the rule of inference $\mathit{I}$ along with other rules.

We fix a sequent calculus for Strong Kleene logic with identity.

\begin{dfn}[$\SK_=$]\label{SequentsSK}
The \Lmc-system $\SK_=$ consists of the following:
	\begin{enumeratei}
		\item the initial sequents $\vphi\Ra\vphi$;
		\item the following inference rules:
		\begin{IEEEeqnarray*}{CCCC}
		\ax{ \fCenter, \vphi }
		\rl{L$\neg$}
			\un{  \neg\vphi, \fCenter }
		\Dis
		& \hspace*{15mm} &
		\axc{ \GD, \vphi }
		\axc{ \vphi, \GD }
		\rl{Cut}
			\bnc{ \GD }
		\Dis
	\\[2mm]
	\ax{ \fCenter }
		\rl{LW}
			\un{  \vphi, \fCenter }
		\Dis
		& \hspace*{15mm} &
		\ax{ \fCenter }
		\rl{RW}
			\un{ \fCenter, \vphi }
		\Dis
	\\[2mm]
		\ax{ \vphi,\fCenter }
		\ax{ \psi,\fCenter }
		\rl{L$\lor$}
			\bn{  \vphi\lor\psi,\fCenter }
		\Dis
		& &
		\ax{ \fCenter,\vphi,\psi }
		\rl{R$\lor$}
			\un{  \fCenter,\vphi\lor\psi }
		\Dis
	\\[2mm]
		\ax{ \vphi,\psi, \fCenter }
		\rl{L$\land$}
			\un{ \vphi\land\psi,\fCenter }
		\Dis
		& &
		\ax{ \fCenter,\vphi }
		\ax{ \fCenter,\psi }
		\rl{R$\land$}
			\bn{ \fCenter,\vphi\land\psi }
		\Dis
	\\[2mm]
		\ax{ \vphi(u), \fCenter }
		\rl{L$\exists$}
			\un{ \exists x\vphi, \fCenter }
			\noLine
		\Dis
		& &
		\ax{ \fCenter, \vphi(t) }
		\rl{R$\exists$}
			\un{ \fCenter, \exists x\vphi }
		\Dis
	\\[2mm]
		\ax{ \vphi(t), \fCenter }
		\rl{L$\forall$}
			\un{ \forall x\vphi, \fCenter }
		\Dis
		& &
		\ax{ \fCenter, \vphi(u) }
		\rl{R$\forall$}
			\un{ \fCenter, \forall x\vphi }
		\Dis
	\\[2mm]
		\axc{ }
		\rl{Ref}
			\un{ \fCenter, t=t }
		\Dis
		& &
		\ax{ \fCenter, \vphi(t) }
		\rl{Repl}
			\un{ \fCenter, s\neq t, \vphi(s) }
		\Dis
	\end{IEEEeqnarray*}
	\end{enumeratei}
\text{Conditions of application: $u$ eigenvariable}
\end{dfn}

Due to the rules Ref (\textit{reflexivity}) and Repl (\textit{replacement}), the system $\SK_=$ derives the sequent $\emptyset\Rightarrow t=s, t\neq s$. Together with the derivability of $t=s, t\neq s \Rightarrow \emptyset$ and the rule Cut, this entails that $\SK_=$ behaves classically on \lnat-formulae:

\begin{obs}
	For $\vphi\in\lnat$, the sequent $\GD, \vphi,\neg\vphi$ is derivable in $\SK_=$.
\end{obs}

Next, we define the theory $\PA[\SK]$:

\begin{dfn}
	The \Lmc-system $\PA[\SK]$ is obtained by expanding $\SK_=$ with the initial sequents of \PA (see e.g.~\cite{takProof}) and the induction rule
	\begin{IEEEeqnarray*}{L}
		\ax{ \vphi(u),\fCenter,\vphi(Su) }
		\rl{{IND}}
			\un{ \vphi(\ovl 0),\fCenter,\vphi(t) }
		\Dis
	\end{IEEEeqnarray*}
for $\vphi(x)\in\Lmc$ and $u$ eigenvariable.
\end{dfn}

By `partial model' we mean a structure $(\Nmc, T, P)$, where: \Nmc is model of arithmetic; $T=(T^+, T^-)$ and $P=(P^+, P^-)$ interpret \T and \P, respectively. The satisfaction relation $(\Nmc, T, P)\svDash{\SK}\vphi$ is defined inductively in the usual way (see, e.g., \cite{kriOutline}). Moreover, we let
	\begin{IEEEeqnarray*}{L}
		(\Nmc, T, P)\svDash{\SK}\GD \text{ iff, if } (\Nmc, T, P)\svDash{\SK}\gamma \text{ for all } \gamma\in\Gamma, \text{ then } (\Nmc, T, P)\svDash{\SK}\delta\text{ for some } \delta\in\Delta.
	\end{IEEEeqnarray*}

We say that $(\Nmc, T, P)$ is a model of a system $S\supset\SK$, written $(\Nmc, T, P)\svDash\SK S$, iff $(\Nmc, T, P)\svDash{\SK}\GD$ whenever $S\vdash\GD$, i.e., whenever $\GD$ is derivable in $S$. A model is \textit{standard} if $\Nmc$ is the standard model \Nbl of arithmetic. Note that models of any system $S\supset\SK$, and in particular models of $\PA[\SK]$, are \textit{consistent}, in the sense that $T^+ \cap T^-=\emptyset$ and $P^+ \cap P^-=\emptyset$. Also, notice that, in an arbitrary model $(\Nmc, T, P)\svDash\SK\PA[\SK]$, if $(\Nmc, T, P)\svDash\SK\vphi\Lra\psi$ and $(\Nmc, T, P)\svDash\SK\neg\vphi\Lra\neg\psi$, then \vphi and $\psi$ have the same semantic value, in the sense that 
	\begin{IEEEeqnarray*}{RCL}
		(\Nmc, T, P)\svDash\SK(\neg)\vphi & \text{ iff } & (\Nmc, T, P)\svDash\SK(\neg)\psi.
	\end{IEEEeqnarray*}
In particular,
	\begin{IEEEeqnarray*}{RCL}
		(\Nmc, T, P)\not\svDash\SK\vphi\lor\neg\vphi & \text{ iff } & (\Nmc, T, P)\not\svDash\SK\psi\lor\neg\psi
	\end{IEEEeqnarray*}
%
We say that $\vphi$ and $\psi$ are \emph{\pask-equivalent} iff $\pask\vdash\vphi\Lra\psi$ and $\pask\vdash\neg\vphi\Lra\neg\psi$.


\section{The theory \tp}

The theory \tp, \emph{Truth and Paradoxicality}, expands $\PA[\SK]$ with principles for truth and paradoxicality. Unless otherwise specified, in what follows we use $s, t$ as variables for codes of closed terms, and we let 
$\vphi, \psi$ range over both formulae and their codes -- context will always make clear which use is intended. So, for example, the sequents
	\begin{IEEEeqnarray*}{L}
		\val s=\val t \Ra \T(\eqn st)
		\\
		\P\vphi\land\P\psi\Ra\P(\vphi\land\psi)
		\\
		\T\forall x\vphi(x) \Ra \forall x\T\vphi(x)
	\end{IEEEeqnarray*}
are shorthand for
	\begin{IEEEeqnarray*}{L}
		\ct u, \ct v, \val u=\val v \Ra \T(\eqn uv)
		\\
		\sent{\and uv}, \P u\land\P v \Ra \P(\and uv)
		\\
		\sent{\all vu}, \T(\all vu) \Ra \forall z\T u(z).
	\end{IEEEeqnarray*}

The following truth principles are essentially the initial sequents of the system \pkf, introduced by Halbach and Horsten \cite{hhAxiomatizing}. The only addition consists of initial sequents ensuring the transparency of \T over the language expanded with \P.

\begin{dfn}[Truth principles]
\hfill
	\begin{IEEEeqnarray*}{RCL+C+L}
		(\mrm{T_1}) & \hspace{3mm} & \val s = \val t \Lra \T(\eqn st) 
				& & \val s \neq \val t \Lra \T(s\ud\neq t)  \hspace*{29mm}
		\\
		(\mrm{T_2}) & \hspace{3mm} & \P\vphi \Lra \T\P\vphi 	
				& & \neg\P\vphi\Lra\T\neg\P\vphi
		\\
		(\mrm{T_3}) & \hspace{3mm} & \T\vphi \Lra \T\T\vphi
				& & \T\neg\vphi \Lra \neg\T\vphi
		\\
		(\mrm{T_4}) & \hspace{3mm} & \T\vphi\land\T\psi \Lra \T(\vphi\land\psi) 
				& & \T\vphi\lor\T\psi \Lra \T(\vphi\lor\psi)
		\\
		(\mrm{T_{5}}) & \hspace{3mm} & \forall x\T\vphi(x) \Lra \T\forall x\vphi(x)
				& & \exists x\T\vphi(x) \Lra \T\exists x\vphi(x)
	\end{IEEEeqnarray*}
\end{dfn}

As per the main points outlined in the introduction, we let $B(x)$ be an \lnat-formula representing the set of \textit{base paradoxical sentences}. The characterization of this set of sentences will of course involve a fairly high degree of arbitrariness and incompleteness -- see \cite[§4.3.1]{casFixed} for a discussion. For the purposes of this short note, however, suffice it to take a sound definition of base paradoxical sentences, including the most common forms of liar sentences, namely: \vphi is a base paradoxical iff it is \pask-equivalent to $\neg\T \vphi$.

To simplify the presentation of the paradoxicality principles, let $\Pi(x):=B(x)\lor B(\ng x)$, and let $A$ be a variable for $\T,\P$.\footnote{
	The formulation of the initial sequents contains redundancies. E.g., since from $(\P_1)$ one can derive $\Pi(\vphi\lor\psi)\Ra\P(\vphi\lor\psi)$, the disjunct $\Pi(\vphi\lor\psi)$ in $(\P_5)$ is redundant in right-to-left direction. Similarly for the remaining principles.
}

\begin{dfn}[Paradoxicality principles]
	\hfill
	\begin{enumerate}[label={$(\P_{\mrm{\arabic*}})$}]
		\item $\Pi(\vphi)\Ra\P\vphi$
		\item $\P\neg A\vphi \Lra \P A\vphi$
		\item $\P\T\vphi \Lra \P\vphi \lor \Pi(\T\vphi)$
		\item $\P(\vphi\land\psi)\Lra(\P\vphi\land\P\psi) \lor (\T\vphi\land\P\psi) \lor (\T\psi\land\P\vphi) \lor \Pi(\vphi\land\psi)$
		\item $\P(\vphi\lor\psi)\Lra(\P\vphi\land\P\psi) \lor (\neg\T\vphi\land\P\psi) \lor (\neg\T\psi\land\P\vphi) \lor \Pi(\vphi\lor\psi)$
		\item $\P\forall x\vphi(x) \Lra \left(\exists x\P \vphi(x) \land \forall y \left(\P\vphi(y)\lor\T\vphi(y)\right)\right) \lor \Pi(\forall x\vphi(x))$
		\item $\P\exists x\vphi(x) \Lra \left(\exists x\P \vphi(x) \land \forall y \left(\P\vphi(y)\lor\neg\T\vphi(y)\right)\right) \lor \Pi(\exists x\vphi(x))$
	\end{enumerate}
\end{dfn}

\begin{dfn}[Interaction principles]
\hfill
	\begin{enumerate}[label={$(\mrm{I}_{\mrm{\arabic*}})$}]
		\item $\T(\vphi\lor\neg\vphi) \Ra \neg\P \vphi$.
	\end{enumerate}
\end{dfn}

\begin{dfn}[\tp]
	The theory \tp is obtained by expanding $\PA[\SK]$ with truth-, paradoxicality-, and interaction-principles.
\end{dfn}

\begin{rmk}
\hfill
	\begin{enumerate}
		\item The disjuncts $\Pi(\xi)$ in $(\P_3)$-$(\P_7)$ capture the \textit{quasi}-compositional account of paradoxicality outlined in the introduction: the paradoxicality of a compound sentence, e.g. $\vphi\land\psi$, is grounded in its components, unless it is a base paradoxical sentence. Note that in $(\P_2)$ the addition of the disjunct $\Pi(\neg A\vphi)$ would be redundant: if $\Pi(\neg A\vphi)$, then $\Pi(A\vphi)$ by definition, hence $\P A\vphi$ by $(\P_1)$.
		\item The contrapositive of $(\mrm{I_1})$, i.e., the initial sequent $\P\vphi \Ra \neg\T(\vphi\lor\neg\vphi)$ is inconsistent over the truth and paradoxicality principles: since $\neg\T(\vphi\lor\neg\vphi)$ is equivalent to $\T(\vphi\land\neg\vphi)$ by $(\mrm{T_3})$, we would obtain that paradoxical sentences are both true and false. For example, let $\lambda$ be a liar sentence such that $\Pi(\lambda)$. Then, from $\emptyset\Ra\P\lambda$ we would get $\emptyset\Ra\neg\T(\lambda\lor\neg\lambda)$ and hence $\emptyset\Ra\T\lambda\land\neg\T\lambda$. However, we will see shortly that a rule encoding an analogous principle is derivable in \tp (Observation~\ref{obs_der}.1).
	\end{enumerate}
\end{rmk}

Before defining a standard model for \tp, let us begin with a few simple observations. First, \tp satisfies the desideratum that \vphi is paradoxical iff its negation $\neg\vphi$ is.

\begin{obs}
	$\tp\vdash\P\vphi\Lra\P\neg\vphi$.
\end{obs}

\begin{proof}
	By formal induction on the buildup of \vphi. If $\vphi$ is atomic, then the claim is just an instance of $(\P_2)$. As an example for a complex sentence, assume by IH that $\P\vphi\Lra\P\neg\vphi$ and $\P\psi\Lra\P\neg\psi$. Then consider the disjunction $\vphi\lor\psi$, for which we want to show
	\begin{IEEEeqnarray*}{L}
		\P(\vphi\lor\psi) \Lra \P(\neg(\vphi\lor\psi)).
	\end{IEEEeqnarray*}
By definition, $\neg(\vphi\lor\psi):=\neg\vphi\land\neg\psi$. One can then reason as follows:
	\begin{IEEEeqnarray*}{RCL+L}
		\P(\vphi\lor\psi) & \Lra & \P(\neg\vphi\land\neg\psi) & \neg(\vphi\lor\psi):=\neg\vphi\land\neg\psi
		\\
		& \Lra & (\P\neg\vphi\land\P\neg\psi) \lor (\neg\T\vphi\land\P\neg\psi) \lor (\neg\T\psi\land\P\neg\vphi) \lor \Pi(\neg\vphi\land\neg\psi) \hspace*{4mm} & (\P_4)
		\\
		& \Lra & (\P\neg\vphi\land\P\neg\psi) \lor (\neg\T\vphi\land\P\neg\psi) \lor (\neg\T\psi\land\P\neg\vphi) \lor \Pi(\vphi\lor\psi) & \Pi(\vphi)\Lra\Pi(\neg\vphi)
		\\
		& \Lra & (\P\vphi\land\P\psi) \lor (\neg\T\vphi\land\P\psi) \lor (\neg\T\psi\land\P\vphi) \lor \Pi(\vphi\lor\psi) & \text{IH}
	\end{IEEEeqnarray*}

The last displayed formula is the consequent of $(\P_5)$.
\end{proof}

\begin{obs}\label{obs_der}
	The following are derivable in \tp:
	\begin{enumerate}
		\item The rule
		\[
			\ax{\fCenter, \P\vphi}
				\un{\T(\vphi\lor\neg\vphi), \fCenter}
			\Dis
		\]
		\item The sequents $\P t\Ra\neg\P\P t$ and $\neg\P t\Ra\neg\P\neg\P t$.
	\end{enumerate}
\end{obs}

\begin{proof}
Both claims follow from $(\mrm{I_1})$. For the rule, we reason as follows:
	\begin{IEEEeqnarray*}{L}
		\axc{ \GD, \P\vphi }
		\ax{ \T(\vphi\lor\neg\vphi), \fCenter, \neg\P\vphi }
			\rl{$\neg$L}
			\dl
			\un{ \P\vphi, \T(\vphi\lor\neg\vphi), \fCenter }
		\bnc{ \T(\vphi\lor\neg\vphi), \GD }
		\Dis
	\end{IEEEeqnarray*}
For the sequents, we have:
	\begin{IEEEeqnarray*}{L}
	\def\fCenter{ \mbox{ $\Ra$ } }
		\ax{ \P t\fCenter\T\P t }
		\dl
		\rl{$\mrm{T_4}$}
			\un{ \P t\fCenter\T(\P t\lor\neg\P t) }
		\ax{ \T(\P t\lor\neg\P t)\fCenter\neg\P\P t }
		\dl
		\bnc{\P t\fCenter\neg\P\P t}
		\Dis
	\end{IEEEeqnarray*}
Similarly for $\neg\P \vphi\Ra\neg\P\neg\P \vphi$.
\end{proof}

\section{Fixed-point semantics}\label{sec:FPS}

We define fixed-point models for \tp.

\begin{dfn}
Let $\mathscr{P}_i(x)$, for $1\leq i\leq 7$, be defined as follows:
	\begin{enumerate}[label={$\mathscr{P}_{\arabic*}(x) :=$}]\setlength{\itemindent}{8mm}
		\item $\sent x\land\Pi(x)$
		\item $\sent x\land\exists s(x=\ud\T s \land \P\val s)$
		\item $\sent x\land\exists s(x=\ng{\ud\T} s \land\P\val s)$
		\item $\sent x \land \exists y\exists z\big[x=\and yz\land \big[ \big(\P x\land\P y\big)\lor \big(\T x\land\P y\big) \lor \big(\T y\land\P x\big) \big]\big]$
		\item $\sent x \land \exists y\exists z\big[x=\dor yz\land \big[ \big(\P x\land\P y\big)\lor \big(\neg\T x\land\P y\big) \lor \big(\neg\T y\land\P x\big) \big]\big]$
		\item $\sent x \land \exists v\exists s\big[x=\all vs\land \exists y\P s(\dot{y}) \land \forall y \big(\P s(y)\lor\T s(y)\big)\big]$
		\item $\sent x \land \exists v\exists s\big[x=\ex vs\land \exists y\P s(\dot{y}) \land \forall y \big(\P s(y)\lor\neg\T s(y)\big)\big]$
	\end{enumerate}
\end{dfn}

It can be observed that $\mathscr{P}(x):=\bigvee_{1\leq i\leq 7}\mathscr{P}_i(x)$ describes the closure conditions of a paradoxicality predicate \P relative to a fixed interpretation of \T. This induces the following operator:
	\[
		\Gamma^+_{\mcr{P},T}(P) := \big\{ \#\vphi\in\omega\mid(\Nbl, T, P)\svDash\SK\mcr{P}(\vphi) \big\}.
	\]
In this sense, the operator $\Gamma^+_{\mcr{P},T}(P)$ yields the set of paradoxical sentences within the structure \ntp{}. In light of the interaction principle $(\mrm{I_1})$, define
	\[
		\Gamma^-_{\mcr{P},T}(P) := \big\{ \#\vphi\in\omega\mid(\Nbl, T, P)\svDash\SK\vphi\lor\neg\vphi\big\}.
	\]
Taken together, $\Gamma^+_{\mcr{P},T}(P)$ and $\Gamma^-_{\mcr{P},T}(P)$ yield:
	\[
		\Gamma_{\mcr{P},T}(P) := \left(\Gamma^+_{\mcr{P},T}(P), \Gamma^-_{\mcr{P},T}(P)\right).
	\]
This is similar to the well-known Kripke Jump, an operator yielding the set of true sentences within a given structure. In the current setting, a suitable Kripke Jump is given by
	\begin{IEEEeqnarray*}{L}
		\Gamma_{\mcr{T},P}(T) := (\Gamma^+_{\mcr{T},P}(T), \Gamma^-_{\mcr{T},P}(T)) := \Big(\{ \#\vphi\in\omega\mid(\Nbl, T, P)\svDash\SK\vphi \}, \ \{\#\vphi\in\omega\mid(\Nbl, T, P)\svDash\SK\neg\vphi\}\Big).
	\end{IEEEeqnarray*}
Combining the two operators together, we get the following:

\begin{dfn}
Let $\Gamma_{\mcr{TP}}:(\mc P(\omega)^2)^2\longrightarrow(\mc P(\omega)^2)^2$ be defined by
 	\begin{IEEEeqnarray*}{RCL}
 		\Gamma_{\mcr{TP}}(T,P) &=& \big(\Gamma_{\mcr{T},P}(T), \Gamma_{\mcr{P},T}(P)\big).
 	\end{IEEEeqnarray*}
\end{dfn}

\begin{dfn}
	For $X_k, Y_k\in\mc P(\omega)$, let $(X_i, Y_i)\leq(X_j, Y_j)$ be defined as $X_i\Sub X_j$ and $Y_i\Sub Y_j$, and $\ang{(X_i, X_j), (Y_i, Y_j)}\leq\ang{(X_l, X_m), (Y_l, Y_m)}$ as $(X_i, X_j)\leq(X_l, X_m) \ \& \ (Y_i, Y_j)\leq(Y_l, Y_m)$.
\end{dfn}

\begin{fact}\label{ft_sk-submodels}
	Let $(T, P)\leq(T', P')$. Then, for all \vphi, if $(\Nbl, T, P)\svDash\SK\vphi$, then $(\Nbl, T', P')\svDash\SK\vphi$. In particular, $(\Nbl, T', P')\svDash\SK\mcr P(\vphi)$ whenever $\ntp{}\svDash\SK\mcr P(\vphi)$
\end{fact}

\begin{proof}
	By induction on the complexity of \vphi.
\end{proof}

Fact~\ref{ft_sk-submodels} yields the following

\begin{lemma}
	The jump operator $\Gamma_{\mcr{TP}}$ is \emph{monotone}, in the sense that 
	\begin{IEEEeqnarray*}{CCL}
		\text{if } (T, P) & \leq & (T', P'), \text{ then }
		\Gamma_{\mcr{TP}}(T,P) \leq \Gamma_{\mcr{TP}}(T',P').
	\end{IEEEeqnarray*}
\end{lemma}

\begin{dfn}\label{df_trans-seq}
	Define a sequence $\ang{T_\alpha, P_\alpha}_{\alpha\in ON}$ as follows:
	\begin{IEEEeqnarray*}{RCLCRCL}
		(T_0, P_0) & := & \ang{(\emptyset, \emptyset),(\emptyset, \emptyset)} &  \\
		(T_{\beta+1}, P_{\beta+1}) & := & \Gamma_{\mcr{TP}}(T_{\beta}, P_{\beta})
		\\
		(T_{\lambda}, P_{\lambda}) & := & \bigcup_{\beta<\lambda}(T_{\beta}, P_{\beta})
	\end{IEEEeqnarray*}
\end{dfn}

Unless otherwise specified, in what follows interpretations indexed by ordinals refer to interpretations in the sequence of Definition~\ref{df_trans-seq} leading to $(\tinf, \pinf)$. 

\begin{rmk}
	Note that, for all $\alpha$, $P_\alpha^-=T^+_{\alpha}\cup T^-_{\alpha}$.
\end{rmk}

By the monotonicity of $\Gamma_{\mcr{TP}}$, we get

\begin{crl}\label{cr_weak-incr}
	The sequence $\ang{T_\alpha, P_\alpha}_{\alpha\in ON}$ of Definition~\ref{df_trans-seq} is such that, for all $\alpha$, $(T_\alpha, P_\alpha)\leq(T_{\alpha+1}, P_{\alpha+1})$.
\end{crl}

\begin{proof}
	Since $(T_0, P_0)$ is the empty interpretation, it is trivially contained in $(T_1, P_1)$. The claim then follows immediately by the monotonicity of $\Gamma_{\mcr{TP}}$ and the definition of $(T_\lambda, P_\lambda)$ for $\lambda$ limit.
\end{proof}

\begin{prop}
	The sequence from Definition~\ref{df_trans-seq} reaches a fixed-point, i.e., a pair $(\tinf,\pinf)$ such that $(\tinf,\pinf) = \Gamma_{\mcr{TP}}(\tinf,\pinf)$. By usual arguments, it can also be shown that $(\tinf,\pinf)$ is the \emph{least} fixed-point of $\Gamma_{\mcr{TP}}$.
\end{prop}

The goal is now to show that $\ntp\infty\svDash\SK\tp$. The bulk of the proof consists in showing the following claims:
	\begin{enumeratei}
		\item \ntp\infty is \emph{consistent}, in the sense that $\tinf^+\cap\tinf^-=\emptyset$ and $\pinf^+\cap\pinf^-=\emptyset$,
		\item \ntp\infty is \emph{sound}, in the sense that $(T^+_\infty\cup T^-_\infty)\cap P^+_\infty=\emptyset$.
	\end{enumeratei}
	
The argument for claim (i) provided below establishes in fact a stronger claim, namely the consistency of every interpretation $(T_\alpha, P_\alpha)$; it can be broken down into the following steps:
	\begin{itemize}
		\item A structure \ntp{} is \textit{inconsistent} -- i.e. either $T^+\cap T^-\neq\emptyset$ or $P^+\cap P^-\neq\emptyset$ -- iff $\ntp{}\svDash\SK\vphi\land\neg\vphi$ for some \vphi.
		\item If \ntp\infty is inconsistent, then there exists a least inconsistent $(T_\alpha, P_\alpha)$. 
		\item Since $T^+_\alpha\cap T^-_\alpha\neq\emptyset$ iff $\ntp\beta\svDash\SK\vphi\land\neg\vphi$ for some $\beta<\alpha$, and since $\ntp\beta\svDash\SK\vphi\land\neg\vphi$ iff $\ntp\beta$ is inconsistent, the least inconsistent $(T_\alpha, P_\alpha)$ is inconsistent in $P_\alpha$ only, that is, $T^+_\alpha\cap T^-_\alpha=\emptyset$ and $P^+_\alpha\cap P^-_\alpha\neq\emptyset$. 
		\item By definition, $\psi\in P_\alpha^+\cap P_\alpha^-$ iff, for $\beta<\alpha$, $\ntp\beta\svDash\SK\mcr P(\psi) \land (\psi\lor\neg\psi)$.
		\item However, it will be shown that if an interpretation \ntp\beta is consistent, then either $\ntp\beta\not\svDash\SK\mcr P(\psi)$ or $\ntp\beta\not\svDash\SK\vphi\lor\neg\vphi$.
	\end{itemize}

We now develop this outline in detail.

\begin{dfn}
	A structure $\ntp{}$, or just an intepretation $(T, P)$, is defined to be 
	\begin{itemize}
		\item \emph{consistent}, if both $T^+\cap T^-=\emptyset$ and $P^+\cap P^-=\emptyset$; \emph{inconsistent} if it is not consistent.
		\item \emph{sound}, if either $(\Nbl, T, P)\not\svDash\SK\vphi\lor\neg\vphi$ or $(\Nbl, T, P)\not\svDash\SK\mcr P(\vphi)$.
	\end{itemize}
\end{dfn}

The next few results are preparatory for the Main Lemma~\ref{le_main}, in which it will be established that consistency entails soundness.


\begin{obs}\label{ob_conistent-models}
	Let $(T, P)$ be consistent. Then, either $(\Nbl, T, P)\not\svDash\SK\vphi$ or $(\Nbl, T, P)\not\svDash\SK\neg\vphi$.
\end{obs}

\begin{proof}
	By induction on the complexity of \vphi.
\end{proof}


\begin{obs}\label{obs_cons}
	If $(T_\alpha, P_\alpha)$ is consistent, then $(T_\gamma, P_\gamma)$ is consistent for every $\gamma<\alpha$.
\end{obs}

\begin{proof}
	By Fact~\ref{ft_sk-submodels} and Corollary~\ref{cr_weak-incr}, inconsistency is preserved upwards.
\end{proof}

\begin{fact}\label{ft_T-InOut}
	For all $\alpha$: If $\ntp\alpha\svDash\SK(\neg)\T\vphi$, then $\ntp\alpha\svDash\SK(\neg)\vphi$.
\end{fact}

\begin{proof}
	Suppose $\ntp\alpha\svDash\SK(\neg)\T\vphi$, which is the case iff $\ntp{\beta}\svDash\SK(\neg)\vphi$ for some $\beta<\alpha$. Since $(T_{\beta}, P_{\beta})\leq(T_\alpha, P_\alpha)$ by Corollary~\ref{cr_weak-incr}, we get $\ntp{\alpha}\svDash\SK(\neg)\vphi$ by Fact~\ref{ft_sk-submodels}.
\end{proof}

Note that, by definition of $\Gamma_{\mcr{TP}}$, the anti-extension $T^-_\alpha$ of $T_\alpha$ can be defined via $T^+_\alpha$ as follows:

\begin{fact}\label{ft_antiext-via-ext}
	For all $\alpha$: $T_\alpha^-=\{\vphi\mid\neg\vphi\in T^+_\alpha\}$. In particular, $\ntp\alpha\svDash\SK\neg\T\vphi\Lra\T\neg\vphi$.
\end{fact}

Since any consistent $\ntp{}$ is a model of \pask, as a corollary we obtain:

\begin{crl}\label{cr_TF}
	If $\Nbl\vDash\Pi(\vphi)$, then $\ntp\alpha\svDash\SK\vphi\Lra\neg\T\vphi$ for any consistent \ntp\alpha.
\end{crl}

\begin{proof}
	Suppose $\Nbl\vDash\Pi(\vphi)$. If $\Nbl\vDash B(\vphi)$, then the claim is obvious. If $\Nbl\vDash B(\neg\vphi)$, then $\ntp\alpha\svDash\SK\neg\vphi\Lra\neg\T\neg\vphi$ iff, by properties of \vphi and definition of double negation, $\ntp\alpha\svDash\SK\vphi\Lra\T\neg\vphi$ iff, by Fact \ref{ft_antiext-via-ext}, $\ntp\alpha\svDash\SK\vphi\Lra\neg\T\vphi$.
\end{proof}

Moreover, since $\Pi(x)\in\lnat$, we have:

\begin{fact}\label{ft_P1alpha}
	$\vphi\in P^+_1$ iff, for all $\alpha$, $\ntp\alpha\svDash\SK\Pi(\vphi)$.
\end{fact}

We thus obtain the following 

%

\begin{lemma}\label{le_liars-tinf}
	Let $(T_\alpha, P_\alpha)$ be consistent. Then $\ntp\alpha\not\svDash\SK\vphi\lor\neg\vphi$ for any $\vphi\in P^+_1$.
\end{lemma}

\begin{proof}
	Let $\vphi\in P^+_1$ and let $\alpha$ be arbitrary. By Fact \ref{ft_P1alpha}, $\ntp\alpha\svDash\SK\Pi(\vphi)$. Towards a contradiction, assume $\ntp\alpha\svDash\SK\vphi\lor\neg\vphi$. This is the case iff
	\begin{IEEEeqnarray*}{L+L}
		\ntp\alpha\svDash\SK\neg\T\vphi\lor\T\vphi, & \text{by } \ntp\alpha\svDash\SK\Pi(\vphi) \text{ and Corollary \ref{cr_TF}},
		\\
		\ntp{\alpha}\svDash\SK\neg\T\vphi\lor\vphi & \text{by Fact \ref{ft_T-InOut}},
		\\
		\ntp{\alpha}\svDash\SK\neg\T\vphi & \text{by } \ntp\alpha\svDash\SK\Pi(\vphi) \text{ and Corollary \ref{cr_TF}},,
		\\
		\ntp{\alpha}\svDash\SK\vphi\land\neg\vphi  & \text{by } \ntp\alpha\svDash\SK\Pi(\vphi), \text{ Fact \ref{ft_T-InOut}, and Corollary \ref{cr_TF}}.
	\end{IEEEeqnarray*}
The last line contradicts the consistency of $\alpha$.
\end{proof}

We can now prove the main lemma, namely, the fact that a structure \ntp\alpha is sound if it is consistent. This amounts to showing that \vphi is undefined in \ntp\alpha whenever $\ntp\alpha\svDash\SK\mcr P(\vphi)$. The claim will be by proven induction on what will be called the \textit{paradoxicality rank} of \vphi, namely, the least $\beta$ such that $\vphi\in P_\beta$.

\begin{dfn}[Paradoxicality rank]
	The \emph{paradoxicality rank}, or $\prn$ for short, of a formula $\vphi\in\pinf^+$ is defined to be least $\alpha$ such that $\vphi\in P^+_\alpha$.
\end{dfn}

\begin{lemma}[Main Lemma]\label{le_main}
	If \ntp\alpha is consistent, then it is sound.
\end{lemma}

\begin{proof}
	Assume $(T_\alpha, P_\alpha)$ is consistent, and let $\ntp\alpha\svDash\SK\mcr P(\vphi)$. We want to show
	\begin{IEEEeqnarray}{L}
		\ntp\alpha\not\svDash\SK\vphi\lor\neg\vphi. \ztag{$*$}
	\end{IEEEeqnarray}	
We reason by induction on $\prn(\vphi)$, which we know must be $\leq\alpha+1$. If $\prn(\vphi) = 1$, then the claim follows from Lemma~\ref{le_liars-tinf}. For $\beta\geq1$, assume the claim holds for sentences of paradoxicality rank $\leq\beta$, and let $\prn(\vphi)=\beta+1$. This entails $\ntp\beta\svDash\SK\mcr P(\vphi)$. Reasoning by cases, show that $\ntp\beta\not\svDash\SK\vphi\lor\neg\vphi$. Here are some examples.

\medskip
\noindent
Suppose $\ntp\beta\svDash\SK\mcrs P2(\vphi)$. Then $\vphi\equiv\T\psi$ and $\ntp{\beta}\svDash\SK\P\psi$. Since $\prn(\psi)\leq\beta$, by IH we have $\ntp\alpha\not\svDash\SK \psi\lor\neg\psi$, hence $\ntp\alpha\not\svDash\SK \T\psi\lor\neg\T\psi$ by Fact~\ref{ft_T-InOut}.

\medskip
\noindent
Similarly if $\ntp\beta\svDash\SK\mcrs P3(\vphi)$. 

\medskip
\noindent
Suppose $\ntp\beta\svDash\SK\mcrs P4(\vphi)$. Then $\vphi\equiv\psi\land\theta$ and one of the following holds
	\begin{IEEEeqnarray*}{LL}
		1. \hspace*{3mm} & \ntp{\beta}\svDash{\SK}\P\psi\land\P\theta
			\\
		2. & \ntp{\beta}\svDash{\SK}\P\psi\land\T\theta
			\\
		3. & \ntp{\beta}\svDash{\SK}\T\psi\land\P\theta.
	\end{IEEEeqnarray*}
If 1., the the claim follows immediately by IH, so suppose 2. Since $\prn(\psi)\leq\beta$, by IH $\ntp{\alpha}\not\svDash\SK\psi\lor\neg\psi$, hence $\ntp{\alpha}\not\svDash\SK\psi\land\theta$. Moreover, from $\ntp{\beta}\svDash{\SK}\T\theta$ we get $\ntp{\alpha}\svDash{\SK}\T\theta$ by upward persistence, and hence $\ntp{\alpha}\svDash{\SK}\theta$ by Fact~\ref{ft_T-InOut}. Since $\ntp\alpha$ is consistent, we get $\ntp{\alpha}\not\svDash{\SK}\neg\theta$, which together with $\ntp{\alpha}\not\svDash{\SK}\neg\psi$ entails $\ntp{\alpha}\not\svDash\SK\neg(\psi\land\theta)$. We conclude $\ntp{\alpha}\not\svDash\SK\psi\land\theta\lor\neg(\psi\land\theta)$. The argument is symmetric for 3.

\medskip
\noindent
Similarly if $\ntp\beta\svDash\SK\mcrs P5(\vphi)\lor\mcrs P6(\vphi)\lor\mcrs P7(\vphi)$.
\end{proof}

A corollary of the Main Lemma is that \ntp{\alpha+1} is consistent whenever \ntp\alpha is. This is because (i) no sentence will be in $T^+_{\alpha+1}\cap T_{\alpha+1}^-$ by consistency of \ntp\alpha and (ii) no sentence will be in $P^+_{\alpha+1}\cap P_{\alpha+1}^-$ by soundness of \ntp\alpha. We then obtain:

\begin{crl}\label{cr_infty-cons}
	\ntp\infty is consistent and (hence) sound.
\end{crl}

\begin{proof}
	By induction on $\alpha$, it can be shown that every \ntp{\alpha} is consistent. The interpretation \ntp0 is vacuously consistent, and the Main Lemma entails that $\ntp{\alpha+1}$ is consistent whenever $\ntp\alpha$ is. For limits $\kappa$, if \ntp\kappa is inconsistent, there must be $\xi, \eta<\kappa$ such that either $T^+_\xi\cap T^-_\eta\neq\emptyset$, or $P^+_\xi\cap P^-_\eta\neq\emptyset$. But then, since the sequence of $\ntp\alpha$ is weakly increasing, there is a $\zeta<\kappa$ such that $T_\xi\cup T_\eta \subseteq T_\zeta$ and $P_\xi\cup P_\eta \subseteq P_\zeta$. This entails that the $\zeta$-th structure \ntp\zeta is inconsistent, contradicting IH.
\end{proof}

We finally want to show that $(\Nbl, \pinf, \tinf)$ is a model of \tp. To this end, we need an auxiliary observation:

\begin{fact}\label{ft_base-liarsP}
	There is no $t$ such that either of the following holds:
	\begin{enumerate}
		\item $\pask\vdash\P t\Lra\neg\T\P t$.
		\item $\pask\vdash\neg\P t\Lra\neg\T\neg\P t$.
		\item $\pask\vdash\T t\Lra\neg\T\T t$.
	\end{enumerate}
\end{fact}

\begin{proof}
	For (1), let $\Mmc:=\langle\Nbl, (\emptyset, \emptyset), (\omega, \emptyset)\rangle$. Then $\Mmc\svDash{\SK}\PA[\SK]$ is such that $\Mmc\svDash{\SK}\P t$, but $\Mmc\not\svDash{\SK}\neg\T\P t$. For (2), let $\Mmc':=\langle\Nbl, (\emptyset, \emptyset), (\emptyset, \omega)\rangle$. Then $\Mmc'\svDash{\SK}\PA[\SK]$ is such that $\Mmc'\svDash{\SK}\neg\P t$, but $\Mmc'\not\svDash{\SK}\neg\T\neg\P t$. For (3), let $\Mmc'':=\langle\Nbl, (\omega, \emptyset), (\emptyset, \emptyset)\rangle$.
\end{proof}

\begin{crl}\label{cr_non-PPt}
	For any $t$, $\{\P t, \neg\P t\}\cap P_\infty^+=\emptyset$.
\end{crl}

\begin{proof}
	Fact \ref{ft_base-liarsP}.1-2 entails $\{\P t,\neg\P t\}\cap P_1^+=\emptyset$. For $\alpha>1$, the claim $\{\P t,\neg\P t\}\cap P_\alpha^+=\emptyset$ is obvious by definition $\Gamma_{\mcr{TP}}$.
\end{proof}

\begin{theorem}\label{th_tp-cons}
	$(\Nbl, \pinf, \tinf)\vDash\tp$, hence \tp is $\omega$-consistent.
\end{theorem}

\begin{proof}
	We show $(\Nbl, \pinf, \tinf)\svDash\SK (\mrm{P_2})$. For $A=\T$, by Fact \ref{ft_base-liarsP}.3, we have that $\T t\in P^+_\infty$ iff $(\Nbl, \pinf, \tinf)\svDash\SK\mcrs P2(\T t)$. This entails $\mrm{val}(t)\in \pinf^+$, hence $\neg\T t\in P^+_\infty$ by $\mcrs P3$. If $\neg\T t\in P^+_\infty$, then $\ntp\infty\svDash\SK\mcrs P3(\neg\T t)\lor\mcrs P1(\neg\T t)$. If $\ntp\infty\svDash\SK\mcrs P3(\neg\T t)$, then $\mrm{val}(t)\in \pinf^+$, hence $\T t\in P^+_\infty$ by $\mcrs P2$. If $\ntp\infty\svDash\SK\mcrs P1(\neg\T t)$, then $\neg\T t\in P^+_1$, which is the case iff $\T t\in P^+_1$. For $A=\P$, \mrm{(P_2)} holds vacuously by Corollary \ref{cr_non-PPt}.
	
	The satisfiability of the remaining  paradoxicality and truth axioms is clear by construction. For example, consider $(\mrm P_4)$. Letting $\Sigma$ abbreviate the consequent of $(\mrm P_4)$, we reason as follows: $\vphi\land\psi\in P_\infty$ iff
	\begin{IEEEeqnarray*}{L+L}
		\vphi\land\psi\in P_{\infty+1} & \text{ by } P_{\infty}=P_{\infty+1}
		\\
		\ntp\infty\svDash\SK\mcrs P1(\vphi\land\psi) \lor \mcr P_4(\vphi\land\psi) & \text{ by definition of } P_{\infty+1}
		\\
		\ntp\infty\svDash\SK\Sigma & \text{ by definition of  \mcrs P1 and \mcrs P4 }
	\end{IEEEeqnarray*}
	
The interaction principle $(\mrm{I_1})$ is also straightforward, since $\tinf^+\cup\tinf^-=\pinf^-$.
\end{proof}
\section{Some additional observations}

We discuss the status of some paradigmatic paradoxical sentences within \tp. Of course, liar sentences ``saying of themselves'' that they are not true are paradoxical by design: these are sentences \vphi \pask-equivalent to $\neg\T\vphi$. Similarly, it is not difficult to see that (Boolean) Curry sentences are paradoxical. For example, given a sentence $\kappa$ \pask-equivalent to $\neg\T(\kappa)\lor0=1$, it is clear that $\kappa$ is also \pask-equivalent to $\neg\T(\kappa)$.

\subsection{McGee}
	Let $f(x,y)$ be a function such that $f(n, \#\vphi)\mapsto\#\T^n\vphi$, where $\T^n\vphi$ abbreviates $n$-many iterations of \T in front of \vphi. The McGee sentence is a sentence $\mu$ \pask-equivalent to
	\[
		\exists x\neg\T\ud f(x, \cor\mu).
	\]
It can be shown that $\mu\notin \pinf$, and hence that its paradoxicality is independent of \tp.

As a preliminary observation, let us note that neither $\mu$ nor its negation are in $P^+_1$:

\begin{fact}\label{ft_mcgee}
	Neither of the following holds:
	\begin{enumerate}
		\item $\pask\vdash\mu\Lra\neg\T\mu$
		\item $\pask\vdash\neg\mu\Lra\neg\T\neg\mu$
	\end{enumerate}
Moreover, there is no $n$ such that one of the following holds
	\begin{enumerate}
	\setcounter{enumi}{2}
		\item $\pask\vdash\T^n\mu\Lra\neg\T\T^n\mu$
		\item $\pask\vdash\neg\T^n\mu\Lra\neg\T\neg\T^n\mu$
	\end{enumerate}
\end{fact}

\begin{proof}
	In each case, it suffices to define a model of \pask not satisfying the equivalence. For instance, for (1) let \Mmc be a standard, consistent structure such that $\mu\in T^+$ and $\T\mu\in T^-$. Then
	\begin{IEEEeqnarray*}{L}
		\Mmc\svDash\SK\neg\T\T\mu, \text{ hence } \Mmc\svDash\SK\mu, \text{ however } \Mmc\svDash\SK\T\mu.
	\end{IEEEeqnarray*}
For (2), let $\Mmc'$ be a standard, consistent structure where everything is true (and hence nothing is false), i.e., $T^+=\omega$ and $T^-=\emptyset$. Then
	\begin{IEEEeqnarray*}{L}
		\Mmc'\svDash\SK\forall x\T\ud f(x, \cor\mu), \text{ hence } \Mmc'\svDash\SK\neg\mu, \text{ however } \Mmc'\not\svDash\SK\neg\T\neg\mu.
	\end{IEEEeqnarray*}
For (3), given an arbitrary $n$, use e.g. the structure where everything is true. For (4), the structure used for (2) works for $n=0$. Otherwise, just let e.g. $\T^{n-1}\mu\in T^-$ and $\neg\T^{n}\mu\in T^+$.
\end{proof}

We now show that $\mu\notin\pinf$. To begin with, the claim $\mu\notin\pinf^-=\tinf^+\cup\tinf^-$ follows from the fact that \ntp\infty is consistent and has a transparent truth predicate. As for the claim $\mu\notin\pinf^+$, assume the contrary towards a contradiction. First, observe that
	\begin{IEEEeqnarray}{L}
		\text{for any } k>0, \prn(\mu)<\prn(\T^k\mu). \label{eq}
	\end{IEEEeqnarray}
This is so because $\T^k\mu\in P_{\alpha+1}^+$ iff $\ntp\alpha\svDash\SK\mathscr{P}_2(\T^{k}\mu)$ iff $\ntp\alpha\svDash\SK \P(\T^{{k}-1}\mu)$.

Now let $\prn(\mu) = \beta+1$. This is the case iff
	\begin{IEEEeqnarray*}{L+L}
		\ntp\beta\svDash\SK\mcrs P7(\mu)
		\\
		\ntp\beta\svDash\SK\exists n \P\neg\T^n\cor\mu \land \forall m ( \P\neg\T^m\cor\mu \lor \neg\T\neg\T^m\cor\mu ) & \text{dfn of \mcrs P7}
		\\
		\ntp\beta\svDash\SK\exists n \P\neg\T^n\cor\mu \land \forall m ( \P\neg\T^m\cor\mu \lor \T\neg\neg\T^m\cor\mu ) & \text{Fact \ref{ft_antiext-via-ext}}
		\\
		\ntp\beta\svDash\SK\exists n \P\neg\T^n\cor\mu \land \forall m ( \P\neg\T^m\cor\mu \lor \T^{m+1}\cor\mu ) & \text{$\neg\neg\vphi:=\vphi$}
		\\
		\forall k\in\Nbl,\ntp\beta\svDash\SK\P\neg\T^k\cor\mu & \mu\notin\tinf^+
	\end{IEEEeqnarray*}
Then, for $k$ arbitrary, from $\ntp\beta\svDash\SK\P(\neg\T^{k}\cor\mu)$ (and Fact \ref{ft_sk-submodels}), we get
	\begin{IEEEeqnarray*}{L}
		\ntp\beta\svDash\SK\mcrs P1(\neg\T^{k}\cor\mu) \lor \mcrs P3(\neg\T^{k}\cor\mu).
	\end{IEEEeqnarray*}
If $\ntp\beta\svDash\SK\mcrs P1(\neg\T^{k}\cor\mu)$, then $\Nbl\vDash\Pi(\neg\T^{k}\cor\mu)$, contradicting Fact \ref{ft_mcgee}. If $\ntp\beta\svDash\SK\mcrs P3(\neg\T^{k}\cor\mu)$, then $\T^{k-1}\mu\in P_\beta$. This entails that $\prn(\T^{k-1}\mu) \leq\beta < \prn(\mu)$, which is impossible by (\ref{eq}).

\begin{rmk}
	By modifying the definition of $B(x)$ along the lines suggested by \cite{casFixed}, one could ensure that $\mu$ turns out paradoxical; we will not do so here.
\end{rmk}

\subsection{Gupta sentence}
	Let $\gamma$ be the Gupta sentence: $\forall x(\T x\lor\neg\T x)$. Observe that $\gamma$ cannot be true or false. To see this, let $u$ be an eigenvariable and $\lambda$ a liar sentence, then:
	
	\begin{IEEEeqnarray*}{L}
	\def\fCenter{ \mbox{ $\Ra$ } }
		\ax{\T u\land\neg\T u \fCenter\emptyset}
			\un{ \exists x(\T x\land\neg\T x) \fCenter\emptyset }
			\un{ \neg\gamma \fCenter\emptyset }
		\ax{\T\lambda\lor\neg\T \lambda \fCenter\emptyset}
			\un{ \forall x(\T x\lor\neg\T x) \fCenter\emptyset }
			\un{ \gamma \fCenter\emptyset }
		\bnc{\gamma\lor\neg\gamma \fCenter\emptyset}
		\Dis
	\end{IEEEeqnarray*}

We now show that it is not paradoxical either. Since clearly $\gamma\notin P_1$, we have $\gamma\in\pinf$ iff
	\begin{IEEEeqnarray*}{L+L}
		\ntp\infty\svDash\SK\mcrs P6(\gamma)
		\\
		\ntp\infty\svDash\SK\exists x \P(\T x\lor\neg\T x) \land \forall y ( \P(\T y\lor\neg\T y) \lor \T(\T y\lor\neg\T y) ) \hspace*{4mm} & \text{dfn of \mcrs P6}
		\\
		\ntp\infty\svDash\SK\exists x \P(\T x\lor\neg\T x) \land \forall y ( \P(\T y\lor\neg\T y) \lor (\T y\lor\neg\T y)) ) & \text{properties of \T in \ntp\infty}
		\\
		\ntp\infty\svDash\SK\exists x \P(\T x\lor\neg\T x) \land \forall y \left( \P(y\lor\neg y) \lor (\T y\lor\neg\T y)) \right) & \text{properties of \P in \ntp\infty}
	\end{IEEEeqnarray*}

However, the second conjunct is falsified by choosing $y:=\tau$ for $\tau$ a truth-teller.

\begin{rmk}
	It is unclear whether the Gupta sentence is or is not paradoxical; we will not discuss this issue here.
\end{rmk}

\subsection{Revenge sentence}
	Let $\rho$ be a revenge sentence which is \pask-equivalent to $\neg\T\cor\rho\lor\P\cor\rho$. It can be shown that $\rho$ is undefined in \ntp\infty, hence independent of the axiom of \tp: since \ntp\infty is sound, we have $\rho\notin\pinf^+$, otherwise $\ntp\infty\svDash\SK\P\rho$, hence $\ntp\infty\svDash\SK\rho$, and therefore $\ntp\infty\svDash\SK\T\rho$, contradicting the soundness of \ntp\infty. Similarly, since \ntp\infty is consistent, $\rho\notin\tinf^+\cup\tinf^-$, hence $\rho\notin\pinf^-$.

\section{Expansions}
	As mentioned in the introduction, there have been proposals attempting to single out properties which demarcate paradoxical sentences from non-paradoxical sentences. Two examples from recent literature are given by Field's \cite{fiePower} and Fujimoto \& Halbach's \cite{fhClassical}. The former develops an \SK-system for truth and a predicate of `strong classicality', which can however also be read as `groundedness' in Kripke's sense \cite{kriOutline} -- cf.~\cite[p.227]{fiePower}. The latter introduces a classical system of truth and determinateness. Both approaches are closely related to each other, and in a sense dual to the approach taken here. 
	
	Both Field and Fujimoto \& Halbach postulate the principle according to which a sentence \vphi is grounded, resp. determinate, iff the sentence ``\vphi is grounded [resp. determinate]'' is grounded, resp. determinate. Of course, a similar principle for paradoxicality, namely $\P t\Lra \P\P t$, is nonsense. Similarly, the principle $\neg\P\neg\P t\Ra\neg\P t$ is not intuitive: e.g., the statement $\neg\P\lambda$ is not paradoxical -- it's just false -- but of course $\lambda$ is paradoxical. Yet, the principle $\neg\P\P t$ may be taken to be intuitive: it may be argued that a statement $\P t$ is never paradoxical, since it is an atomic statement which is not a base paradoxical sentence.
	
	Without arguing for or against this view, we sketch how to modify the construction from §\ref{sec:FPS} so to obtain a model for the theory $\tp^+$, obtained by expanding \tp with the initial sequents $\emptyset\Ra\neg\P\P t$. Recalling that there is no liar-sentence of the form $\P t$ or $\neg\P t$ (Fact \ref{ft_base-liarsP}), the idea to construct a model for $\tp^+$ with essentially the same strategy as the one adopted for constructing the model for \tp.

Let $P_0^{*-}:=\{\P\vphi\mid\vphi\in\Lmc\}$. Define the Jump $\Gamma^*_{\mcr T \mcr P}$ just as the Jump $\Gamma_{\mcr T \mcr P}$ above, except that the set $P_0^{*-}$ is kept along the way in the anti-extension of $P$. That is:
	\begin{IEEEeqnarray*}{L}
		\Gamma^{*-}_{\mcr P, T}(P):=\{\vphi \mid \ntp{}\svDash\SK\vphi\lor\neg\vphi\}\cup P_0^{*-}.
	\end{IEEEeqnarray*}

Since $\Gamma^*_{\mcr{TP}}$ is clearly monotone, it can be shown that it has a minimal fixed-point $\ntp{\infty}^*$.

\begin{dfn}\label{df_trans-seq2}
	For each ordinal $\alpha$, define $P_\alpha$ and $T_\alpha$ as follows
	\begin{IEEEeqnarray*}{RCLCRCL}
		(T_0, P_0)^* & := & \ang{(\emptyset, \emptyset),(\emptyset, \{\P\vphi\mid\vphi\in\Lmc\})} &  \\
		(T_{\xi+1}, P_{\xi+1})^* & := & \Gamma^*_{\mcr{T}, \mcr{P}}(T_{\xi}, P_{\xi})^*
		\\
		(T_{\lambda}, P_{\lambda})^* & := & \bigcup_{\xi<\lambda}(T_{\xi}, P_{\xi})^*
	\end{IEEEeqnarray*}
\end{dfn}

\begin{lemma}\label{le_main2}
	For all $\alpha$, $\ntp\alpha^*\not\svDash\SK\mcr P(\P t) \lor \mcr P(\neg\P t)$.
\end{lemma}

\begin{proof}
	By definition, $\ntp\alpha^*\svDash\SK\mcr P(\P t)$ entails $\ntp\alpha^*\svDash\SK\mcrs P1(\P t)$, hence $\Nbl\vDash\Pi(\P t)$, which is impossible by Fact \ref{ft_base-liarsP}. Similarly for $\mcr P(\neg\P t)$.
\end{proof}

\begin{crl}
	For all $\alpha$ and all $t$, $\ntp\alpha^*\not\svDash\SK\P\P t\lor\P\neg\P t$. In particular, $P_0^{*-}\cap P_\infty^{*+}=\emptyset$.
\end{crl}

The next lemma is the analogous to Lemma \ref{le_main}: it is needed to show that consistency is preserved throughout the sequence leading to the least fixed-point of $\Gamma^*_{\mcr{TP}}$. To this end, we adapt the definition of soundness: $\ntp\alpha^*$ is defined to be \emph{sound\,$^*$} iff $\ntp\alpha^*\not\svDash\SK\mcr P(\vphi)$ whenever either $\vphi\in P^{*-}_0$ or $\ntp\alpha^*\svDash\SK\vphi\lor\neg\vphi$.

\begin{lemma}\label{le_minor2}
	If $\ntp\alpha^*$ is consistent, then it is sound\,$^*$.
\end{lemma}

\begin{proof}[Proof Sketch]
	Assume $(T_\alpha, P_\alpha)^*$ is consistent, and let $\ntp\alpha^*\svDash\SK\mcr P(\vphi)$. By Lemma \ref{le_main2}, $\vphi\notin P^{*-}_0$. It then suffices to show
	\begin{IEEEeqnarray}{L}
		\ntp\alpha\not\svDash\SK\vphi\lor\neg\vphi. \ztag{$*$}
	\end{IEEEeqnarray}
The argument is by induction on $\prn(\vphi)$, and it follows the blueprint of the argument for Lemma \ref{le_main}. 
\end{proof}

\begin{crl}
	$\ntp\infty^*$ is consistent and (hence) sound\,$^*$.
\end{crl}

\begin{proof}[Proof Sketch]
	The above Lemma entails that $\ntp{\alpha+1}^*$ is consistent whenever $\ntp\alpha^*$ is. The argument is simlar to that of Corollary \ref{cr_infty-cons}.
\end{proof}

\begin{theorem}
	$(\Nbl, \pinf, \tinf)^*\vDash\tp^+$.
\end{theorem}

\begin{proof}[Proof Sketch]
	For $(\mrm P_2)$, note that the argument is the same as that for Theorem \ref{th_tp-cons}. In particular, Fact \ref{ft_base-liarsP} entails that $\{\P t, \neg\P t\}\cap P_\infty^{*+}=\emptyset$. Similarly for $(\mrm P_3)$-$(\mrm P_7)$. For the interaction principle $(\mrm{I_1})$, just observe that $\tinf^{*+}\cup\tinf^{*-}\subseteq\pinf^{*-}$. Finally, it is clear by construction that $\ntp\infty^*\svDash\SK\neg\P\P t$.
\end{proof}

\vspace*{20mm}
\bibliography{bib}

\providecommand{\noopsort}[1]{}
\begin{thebibliography}{Tak87}

\bibitem[Bac15]{bacCan}
Andrew Bacon.
\newblock Can the classical logician avoid the revenge paradoxes?
\newblock {\em Philosophical Review}, 124(3):299--352, 2015.

\bibitem[Cas21]{casFixed}
Luca Castaldo.
\newblock Fixed-point models for paradoxical predicates.
\newblock {\em The Australasian Journal of Logic}, 18(7):688--723, 2021.

\bibitem[CN24]{cnClassical}
Luca Castaldo and Carlo Nicolai.
\newblock On classical determinate truth.
\newblock {\em arXiv preprint arXiv:2409.04316}, 2024.

\bibitem[CR25]{crAsymmetries}
Luca Castaldo and Lucas Rosenblatt.
\newblock Asymmetries in naive paradoxicality: Classical vs. nonclassical
  logic.
\newblock {\em Manuscript}, 2025.

\bibitem[Fef08]{fefAxioms}
Solomon Feferman.
\newblock Axioms for determinateness and truth.
\newblock {\em The Review of Symbolic Logic}, 1(2):204--217, 2008.

\bibitem[FH24]{fhClassical}
Kentaro Fujimoto and Volker Halbach.
\newblock Classical determinate truth {I}.
\newblock {\em The Journal of Symbolic Logic}, 89(1):218--261, 2024.

\bibitem[Fie22]{fiePower}
Hartry Field.
\newblock The power of naive truth.
\newblock {\em The Review of Symbolic Logic}, 15(1):225--258, 2022.

\bibitem[HH06]{hhAxiomatizing}
Volker Halbach and Leon Horsten.
\newblock Axiomatizing {K}ripke's theory of truth.
\newblock {\em The Journal of Symbolic Logic}, 71(2):677--712, 2006.

\bibitem[IS23]{isReasoning}
Luca Incurvati and Julian~J. Schl\"odel.
\newblock {\em Reasoning with Attitude: Foundations and Applications of
  Inferential Expressivism}.
\newblock Oxford University Press, 2023.

\bibitem[Kri75]{kriOutline}
Saul Kripke.
\newblock Outline of a theory of truth.
\newblock {\em The Journal of Philosophy}, 72(19):690--716, 1975.

\bibitem[Lei05]{leiDepends}
Hannes Leitgeb.
\newblock What truth depends on.
\newblock {\em Journal of Philosophical Logic}, 34(2):155--192, 2005.

\bibitem[MR20]{mrGeneralized}
Julien Murzi and Lorenzo Rossi.
\newblock Generalized revenge.
\newblock {\em Australasian Journal of Philosophy}, 98(1):153--177, 2020.

\bibitem[Pic19]{picAlethic}
Lavinia Picollo.
\newblock Alethic reference.
\newblock {\em Journal of Philosophical Logic}, pages 1--22, 2019.

\bibitem[Pic20]{picReferenceAnd}
Lavinia Picollo.
\newblock Reference and truth.
\newblock {\em Journal of Philosophical Logic}, 49(3):439--474, 2020.

\bibitem[Rei85]{reiRemarks85}
William Reinhardt.
\newblock Remarks on significance and meaningful applicability.
\newblock {\em Mathematical Logic and Formal Systems}, 94:227, 1985.

\bibitem[Rei86]{reiRemarks}
William Reinhardt.
\newblock Some remarks on extending and interpreting theories with a partial
  predicate for truth.
\newblock {\em Journal of Philosophical Logic}, 15(2):219--251, 1986.

\bibitem[RG02]{rgConceptions}
Lucas Rosenblatt and Camila Gallovich.
\newblock Conceptions of paradoxicality.
\newblock In Lorenzo Rossi, editor, {\em The Liar Paradox}. Cambridge
  University Press, fortcoming, 202?

\bibitem[RG22]{rgParadoxicality}
Lucas Rosenblatt and Camila Gallovich.
\newblock Paradoxicality in {K}ripke’s theory of truth.
\newblock {\em Synthese}, 200(2):71, 2022.

\bibitem[Ros22]{rosShould}
Lucas Rosenblatt.
\newblock Should the non-classical logician be embarrassed?
\newblock {\em Philosophy and Phenomenological Research}, 104(2):388--407,
  2022.

\bibitem[Ros23]{rosParadoxicality}
Lucas Rosenblatt.
\newblock Paradoxicality without paradox.
\newblock {\em Erkenntnis}, 88(3):1347--1366, 2023.

\bibitem[Sch14]{schAxioms}
Thomas Schindler.
\newblock Axioms for grounded truth.
\newblock {\em Review of Symbolic Logic}, (1):73--83, 2014.

\bibitem[Tak87]{takProof}
Gaisi Takeuti.
\newblock {\em Proof Theory}, volume~81 of {\em Studies in Logic and the
  Foundations of Mathematics}.
\newblock Elsevier Science Publishers, second edition, 1987.

\bibitem[TS00]{tsBasic}
Anne~Sjerp Troelstra and Helmut Schwichtenberg.
\newblock {\em Basic proof theory}.
\newblock Number~43. Cambridge University Press, 2000.

\end{thebibliography}
\bibliographystyle{alpha}
\end{document}